\documentclass[11pt]{amsart}
\usepackage{amsbsy,amssymb,amscd,amsfonts,latexsym,amstext,delarray,
amsmath,graphicx,slashed} 
\usepackage[small,nohug,heads=littlevee]{diagrams}
\diagramstyle[labelstyle=\scriptstyle]
\usepackage[all,cmtip]{xy}
\usepackage{fullpage}
\usepackage[pdftex]{color}

\newtheorem{thm}{Theorem}[section]
\newtheorem{prop}[thm]{Proposition}
\newtheorem{cor}[thm]{Corollary}
\newtheorem{lem}[thm]{Lemma}

\newtheorem{rem}[thm]{Remark}

\numberwithin{equation}{section}

\def\bL{{\mathbb L}}

\def\A{{\mathbb A}}
\def\C{{\mathbb C}}

\renewcommand{\P}{{\mathbb P}}
\def\Q{{\mathbb Q}}
\def\Z{{\mathbb Z}}
\def\R{{\mathbb R}}
\def\K{{\mathbb K}}

\def\cA{{\mathcal A}}
\def\cB{{\mathcal B}}

\def\cD{{\mathcal D}}

\def\cH{{\mathcal H}}

\def\cM{{\mathcal M}}
\def\cN{{\mathcal N}}

\def\cQ{{\mathcal Q}}
\def\cR{{\mathcal R}}

\def\cW{{\mathcal W}}

\def\m{{\mathfrak m}}

\title{Position Space Feynman quadrics and their motives  }
\author{\"Ozg\"ur Ceyhan}
\address{Facult\'e des Sciences, de la Technologie et de la Communication, University of Luxembourg, 6 rue Richard 
Coudenhove-Kalergi, L-1359 Luxembourg}
\email{ozgur.ceyhan@uni.lu}

\date{\today}
 
\begin{document}
\maketitle

{\it To Yuri Ivanovitch Manin, on the occasion of his birthday, with admiration}

\begin{abstract}
  In this note, we introduce and study  position space Feynman quadrics that are
    the loci of divergences of the position space Feynman integrals
  for Euclidean massless scalar quantum field theories. We prove that
  the   Feynman quadrics define objects in the category of
  mixed Tate motives  for complete graphs with a bound
    on the number of vertices.  This result  shows a
  strong contrast with the graph hypersurfaces approach
   which produces also non-mixed  Tate
    examples.
\end{abstract}

\section{Introduction}

\subsection{The problem: From Feynman integrals to motives and periods} 
A {\it Feynman graph} is a finite  1-dimensional connected CW-complex. For a
  given Feynman graph $\Gamma$, we denote the sets of vertices and edges respectively by 
$Ver(\Gamma)$ and $Edg(\Gamma)$, and  the boundary map by
 $\partial_\Gamma: Edg(\Gamma) \to Ver(\Gamma)$.

In quantum field theory (QFT), the protocol called {\it Feynman rules} associates an integral to 
each Feynman graph $\Gamma$ of the QFT. In position  space setting  Feynman rules work 
as follows: Let $X$  be the spacetime manifold and $d {\bf x}$ be a volume form on X. 
\begin{itemize}
\item First, take the space  $X^{Ver(\Gamma)}:=\{f:  Ver(\Gamma) \to X \text{ where } v \mapsto 
{\bf x}_v \}$ of possibly degenerate configurations of vertices in $X$;

\item Next, attach the {\it propagator} $G_\R ({\bf x}^i , {\bf x}^j)
  $ of the QFT (which is simply the Green's function
  for the Laplacian) to each edge $e$ with $\partial_\Gamma (e) =
  (ij)$.
\end{itemize}
The {\it (unregularized) Feynman integral} associated to $\Gamma$ is defined as the integral 
\begin{equation}  \label{eqn_feynman_int}
^0 W_\Gamma := \int_{X^{Ver(\Gamma)}} \omega_\Gamma
\end{equation}
of the differential form
\begin{equation}  \label{eqn_feynman_form}
\omega_\Gamma := \prod_{ {e \in Edg(\Gamma)} \atop  {\partial_\Gamma (e)  = (ij)}}
G_\R ({\bf x}^i , {\bf x}^j)  \bigwedge_{v \in Ver(\Gamma)} \ d {\bf x}_v.
\end{equation}
These integrals are generally divergent (see for instance
\cite{IZ}). In order to extract meaningful quantities,
one takes care of the divergencies of the integrand
  \eqref{eqn_feynman_form} by using a {\it regularization} procedure
  which carefully introduces {\it counterterms} to produce the {\it
  regularized values} $W_\Gamma$ of $^0 W_\Gamma$.

Kontsevich was the first who suggested that one should consider
appropriate regularizations of these integrals, as
  well as their residues, as periods in the sense of \cite{KZ}. 
  In this perspective, one expects to apprehend the nature of
the numbers arising from Feynman integrals by examining the 
locus $Z_\Gamma := \{ \omega_\Gamma  = \infty \}$ of
the divergence 
of  the integrand \eqref{eqn_feynman_form} as an algebraic variety in an appropriate complexification 
of the configuration space $X^{Ver(\Gamma)}$.

The first guess was that the motives of  $Z_\Gamma$ are mixed Tate motives for Euclidean 
massless scalar QFTs where the spacetime $X$ is $\R^{2d}$ ($d>1$) and the propagator is
\begin{equation} \label{eqn_euclidean}  
G_\R ({\bf x} , {\bf y}) :=   G_\R ({\bf x} - {\bf y}) =  || {\bf x} - {\bf y}  ||^{2-2d}.
\end{equation}

 \smallskip
 \noindent 
 {\bf Main Problem:}  {\it  Confirm whether  $Z_\Gamma$ define objects in the category of mixed Tate motives.}
 
 \smallskip

 Attacking the main problem is a plausible strategy
in explaining the presence of multiple
 $\zeta$-values in Feynman integral computations (see \cite{BrK}), as it would follow
 from the result of F. Brown \cite{Brown} combined with confirmation of 
 $Z_\Gamma$ are in the category of mixed Tate motives.

\subsection{Main results: Feynman quadrics and  their motives}  
The main focus of this work to examine the geometry of the {\it Feynman quadric} 
\begin{equation} \label{eqn_graph_hypersurface}
Z_\Gamma := \left\{  ({\bf x}^i =(z_1^i,\dots,z_{d}^i,w_1^i,\dots,w_{d}^i) : i \in Ver(\Gamma) ) \in (\A^d \times \A^d)^{Ver(\Gamma)} \mid 
 \prod_{ {e \in Edg(\Gamma)} \atop {\partial e = (ij)}} G ({\bf x}^i - {\bf x}^j) = \infty \right \}  
\end{equation}
where
\begin{equation} \label{eqn_propagator}
 G ({\bf x)} =  \frac{1}{(z_1 \cdot w_1 + \cdots + z_{d} \cdot w_{d} )^{d-1}}. 
\end{equation}
This is a suitable setup to attack the main problem since the fixed point set of the real structure
 $$
 c: 
 \begin{array}{ccc}
 (\A^{2d}(\C) )^{Ver(\Gamma)} & \to & (\A^{2d}(\C) )^{Ver(\Gamma)}  \\
(z_1^i,\dots,z_{d}^i,w_1^i,\dots,w_{d}^i) &\mapsto & (\bar w_1^i,\dots, \bar w_{d}^i , \bar z_1^i,\dots, \bar  z_{d}^i)
\end{array}
\ \ \ \text{for all} \ i \in  Ver(\Gamma),
$$
can be identified with the configuration space  $(\A^{2d}(\R) )^{Ver(\Gamma)}$ in the spacetime  $\R^{2d}$ and the restriction 
of \eqref{eqn_propagator} to the fixed point set of $c$  is 
the Euclidean massless scalar propagator $G_\R$ in \eqref{eqn_euclidean}.


 \smallskip
 \noindent 
 {\bf Main Theorem:}  {\it 
Let $\kappa_{n}$ denote the complete graph with $|Ver(\kappa_{n})| = n$, and let  $n \leq d+1$. Then, the Feynman quadric  
$ Z_{\kappa_{n}}$ defines an object in the category of mixed Tate motives. } 
 
 \smallskip

Considering only  the  complete graphs with limited number of vertices  may seem  very restrictive at the first glance,
but we can still consider the Feynman integrals of  arbitrarily  complicated graphs in higher dimensional spacetimes: 
For a given $\Gamma$
with $|Ver(\Gamma)| = n$,  since the graph $\Gamma$ is a subgraph of $\kappa_n$,   the integrand  
\eqref{eqn_feynman_form}  can be thought as an algebraic form
\begin{equation} \label{eqn_integral_complete} 
\prod_{ {e \in Edg(\Gamma)} \atop  {\partial_\Gamma (e)  = (ij)}}
G ({\bf x}^i -  {\bf x}^j)  \bigwedge_{v \in Ver(\Gamma)} \ d {\bf x}_v \in \Omega^*((\A^{2d})^{n} \setminus Z_{\kappa_n}) 
\ \ \text{with} \ \ d {\bf x}_v :=  dz_1^v \wedge \cdots \wedge dz_d^v \wedge    dw_1^v \wedge \cdots \wedge dw_d^v
\end{equation}
 whenever $ d \geq n-1$.  In other words, after an appropriate regularization\footnote{One needs here a regularization
 protocol that preserves  the divergence loci, such as the algebraic regularization \cite{CeyMar3, MarNi} or
 analytic regularization, see for instance \cite{Speer}.}, 
 the residues and
 the regularized values of these Feynman integrals can be given as periods of mixed Tate motives.

 This technical result actually  implies  a bit more: 
The valid theories describing the same phenomenon must produce the same
outcomes. 
This thesis implies that  both momentum and position space formulations of 
Feynman integrals  should produce the  same results  for any valid regularization 
techniques. The resulting nature of the motives should 
be  independent from its formulation.   However,
our main results above indicates otherwise.

\subsection{Graph hypersurfaces vs Feynman quadrics}

In its  original form in \cite{K97}, the main problem  refers to the motives of the {\it graph hypersurfaces}
$$
X_\Gamma := \left\{ \alpha := (\alpha_e : e \in Edg(\Gamma))   \mid
 \Psi_\Gamma(\alpha) : = \sum_{ {\text{spanning trees}} \atop  {T \subset \Gamma}  } 
 \prod_{e \not\in Edg(T) } \alpha_e = 0 \right\} \subset \A^{Edg(\Gamma)}
$$
that are the loci of the divergences for the {\it parametric Feynman
  integrals} in momentum space.  
  The attempts to prove that the graph
hypersurfaces are mixed Tate motives failed already at early stages
\cite{BB},  and the recently found counter examples in
  \cite{BD, BS, D} quashed all remaining hope for the mid-dimensional
  (co)homology, the part responsible of Feynman integrals, being mixed
  Tate. (See \cite{Mar} and the references therein).

Our main theorem here, and the recent progresses \cite{BD, BS, D} on graph
hypersurfaces,  indicate that the motivic nature of
Feynman integrals depends on the momentum/position space formulations
and the regularization methods. The periods of non-mixed Tate
motives are expected to be complicated in general 
(see, for instance  \S 1.2 in \cite{BS1}). However, the periods of the Feynman quadrics should be
given in term of   multiple polylogarithms according to  Conjecture 1.9 in \cite{Go}.

The key observation here is  that the definition of $X_\Gamma$ is
independent of the dimension of the spacetime. This fact allows us
choose the dimension of the Feynman quadrics high enough so that they
would provide mixed Tate motives as stated in our main theorem.

\subsection{Organization of the paper}

In Section \ref{sec_propagator}, we examine the simplest possible case, that is the singular quadric 
$\{ z_1 \cdot w_1 + \cdots + z_{d} \cdot w_{d}  =0\} \subset \A^{2d}$ associated to
the propagator  \eqref{eqn_propagator}. In this test example, we illustrate the main techniques;
we compute the class of this quadric in Grothendieck ring of varieties as well as the associated mixed 
motive in the Voevodsky's category. Our computations show that  this quadric defines a mixed 
Tate motive. 
After resolving the  singularity of the quadric, we show that 
its motive is a mixed Tate motive. 

In Section \ref{sec_quadric}, we introduce the Feynman quadric $  Z_\Gamma$. We then  
examine the Feynman quadrics $  Z_{\kappa_n }$ associated to the 
complete graphs $\kappa_n$ in details. We give a stratification of   $  Z_{\kappa_n }$  by using 
appropriate projections which   we illustrate on a simple case in \S \ref{sec_propagator}. 
While most of the strata are given by trivial fibrations, one of the strata turn out be highly non-trivial, a {\it space
of hyperplane arrangements in  almost general positions}.

In Section \ref{sec_motive}, we prove that the each of these strata of  $  Z_{\kappa_n }$ define
mixed Tate motives, and concluded  that,  the Feynman quadric itself  defines  also an element in the subcategory
of mixed Tate motives when $n \leq d+1$.

In Section \ref{sec_corollaries}, we discuss the similarities and the variations between the Feynman quadrics
and the alternative geometries arising from the Feynman integrals. We briefly discuss the possible reasons 
behind the contrasts between the graph hypersurfaces and the Feynman quadrics. We also discuss the 
possible generalizations to the cases of more general graphs and  the metrics with different signature. 
Finally, we note other alternatives geometries arising from position space Feynman integrals. 




\tableofcontents

\subsection{Notation and Conventions} 
We denote the category of schemes  of finite type over $\Bbbk$ by  $\mathfrak{Sch}_{\Bbbk}$, and
the Grothendieck ring of varieties by  $K_0(Var)$.
We  use the same notation with Voevodsky \cite{Vo}  to  denote the motivic 
categories over $\Bbbk$, such as $DM^{eff}_{gm}(\Bbbk)$, $DM^{-}_{gm}(\Bbbk)$, 
$DM^{eff}_{gm}(\Bbbk) \otimes \Q$ etc., and  denote the motivic functor (resp. with compact support) 
by $\m: \mathfrak{Sch}_{\Bbbk} \to DM^{eff}_{gm}(\Bbbk)$ 
(resp.   $ \m^c:  \mathfrak{Sch}_{\Bbbk} \to  DM^{eff}_{gm}(\Bbbk)$).

\subsubsection{ The ambient space}
We fix the dimension $d \geq 1$ and consider the product $\P^d \times \P^d$ of
projective spaces with     homogeneous  coordinates 
$([{\bf z} ] ; [{\bf w}]) := ( [ z_0, \dots ,z_d ] ;  [ w_0, \dots ,w_d]  )$. 
We set   $\A^d \times \A^d \hookrightarrow \P^d \times \P^d: ( {\bf z} ; {\bf w}) := 
(z_1 , \dots , z_d ;  w_1 , \dots , w_d)
\mapsto  
( [1: z_1: \dots :z_d] ; [1: w_1:\dots :w_d] ) $ be the embedding of the affine part
 $\A^d \times \A^d$.

\section{Simplest example: A quadric associated to the massless  propagator }
\label{sec_propagator}

In this section, we associate a {\it simple Feynman quadric} to the propagator of the Euclidean 
massless scalar QFTs. While  the simple Feynman quadric is defined as
the locus of divergence of the propagator, its geometric and motivic properties are slightly different 
than its relatives in the literature  \cite{BEK, CeyMar1, CeyMar2,CeyMar3, Mar, MarNi}.

We study the geometry of this quadric  in this section. We  compute its class in the Grothendieck 
ring of varieties and its motive in the Voevodsky's
category of mixed motives. Our computations  conclude  that the simple Feynman quadric defines 
a mixed Tate motive  in Voevodsky's category.  We then resolve the singularity of this  quadric and 
show that its  motive is also mixed Tate.

\subsection{Euclidean scalar massless propagator and the simple Feynman  quadric} Let $d >1 $.
Our centre of interest  will be the   quadric $\cQ \subset  \A^d \times \A^d$ defined by
\begin{equation} 
\label{eqn_quadric}
q( {\bf z} , {\bf w}):= {\bf z} \cdot {\bf w} = z_1 w_1 + \cdots + z_d w_d = 0.
\end{equation}
We call it {\it simple Feynman quadric}. 
We denote the  Zariski closure of $\cQ$   in $\P^d \times \P^d$   by $\widehat{\cQ} $.

\begin{rem}
\label{rem_orthogonal}
The   simple Feynman  quadric   parameterizes the pairs of orthogonal vectors in the vector
space $\K^d$, i.e.,   ${\bf z} , {\bf w} \in \K^d$  and ${\bf z} \perp {\bf w} $. The complement of the
locus $\{ {\bf z} =0 \} \cup \{ {\bf w =0} \}$ is closely related to the
Stiefel  manifolds and can be viewed as a homogeneous space. Such a setup would provide an
elegant way to examine motive of $\cQ$ in terms of the motives of the classical groups. 
A related approach can be found in \cite{MarNi}. 
However, we   follow a more pedestrian approach in this paper. 
\end{rem}

\subsubsection{From   Feynman quadric to the propagator}
 Let  $c: \A^{2d} (\C) \to  \A^{2d} (\C) $   be the real  structure 
\begin{equation} 
\label{eqn_real}
 ( {\bf z} ;  {\bf w} ) \mapsto ( {\bf \bar{w}} ; {\bf \bar{z}}   ) := ( \bar{w}_1 , \dots , \bar{w}_d   ;  
 \bar{z}_1, \dots , \bar{z}_d  ).
\end{equation}

\begin{lem}
The restriction of the rational function $q^{1-d}$ to the real locus 
${\bf Fix}(c)$ of   \eqref{eqn_real}  gives the propagator \eqref{eqn_euclidean} of the Euclidean 
massless  scalar  QFT in $\A^{2d}(\R)$.
\end{lem}

\begin{proof}
The fixed point locus ${\bf Fix}(c)$ of \eqref{eqn_real} is  
$\A^d(\C) = \{(  {\bf z} ,  {\bf w}  )  \mid w_i = \bar{z}_i \  \forall \  i=1,\dots,d \}$, that is  $\A^{2d}(\R)$
as  a smooth manifold. 
 Therefore, the restriction of the function $q^{1-d}$ 
to real locus becomes
$$
\frac{1}{(q(  {\bf z}  ,  {\bf \bar{z}}  ))^{d-1}} =  \frac{1}{(|z_1|^2 + \cdots + |z_d|^2 )^{d-1}}    
= \frac{1}{\| {\bf z} \|^{2d-2}}     
$$
that is exactly the propagator for the scalar massless QFT in $\A^{2d}(\R)$ (see, for instance \S 7 in \cite{GJ}). 
\end{proof}

\subsection{The simple Feynman quadric  and its motive  }



\label{sec_projection}
Consider the projection $ \pi: \A^d \times \A^d  \to  \A^d: (  {\bf z}  ,  {\bf w}  ) \mapsto   {\bf w}  $,
and its  restriction 
\begin{equation} \label{eqn_fibration}
 \pi:   {\cQ}  \to \A^d
\end{equation} 
to  our  quadric $ {\cQ}$.

\begin{lem} \label{lem_projective_fibre}
The  fibres of the morphism   $  {\cQ}  \to \A^d$  are
 \begin{eqnarray*} \label{eqn_projective_fibre}
\pi^{-1}( {\bf w}  ) \cong 
\left\{
	\begin{array}{lcl}	
		\A^{d} = \{ ( z_1,\dots ,z_d )    \}
					&	\text{if} &   {\bf w}  = (0, \dots, 0 ), \\
		\A^{d-1} = \{ ( z_1,\dots ,z_d )   \mid  z_1 w_1 + \cdots + z_d w_d  =0 \}
					&	\text{if} &  {\bf w}  \ne (0, \dots, 0 ). 		
	\end{array}	
 \right.  
\end{eqnarray*}
\end{lem}

\begin{proof}
The statement directly follows from the defining equation \eqref{eqn_quadric} of the
quadric $\cQ$. 
\end{proof}

The types of the fibres in Lemma \ref{lem_projective_fibre} of the projection $ \pi$ induce the following
decomposition of the  base $\A^d$: Let
\begin{eqnarray*}  
\begin{array}{ccl} \label{eqn_cellular}
  	  S_0 &=& \{ (0, \dots, 0 )  \} \subset \A^d = \{ {\bf w} \} \\ 
	  S_1   &= &\A^d \setminus   S_0   = \A^d \setminus \{0\}, \ 
\end{array}
\end{eqnarray*}
and, consider the  fibrations
\begin{eqnarray} \label{eqn_cell_fibration}
\begin{diagram}
		  F_i	&\rInto	&   U_i 	& = 	&   \pi^{-1}(S_{i}) \\
		\dTo_{ \pi} 	&		& \dTo_{ \pi} 	&	&	\\
		\text{pt}		&\rInto	&   S_i 	& 	& 	\
\end{diagram}
\end{eqnarray}
with  fibres $  F_i$ given in Lemma \ref{lem_projective_fibre}, i.e., 
$  F_0 = \A^{d}$ and $  F_1   =  \A^{d-1}$.

\begin{lem} \label{lem_trivial0}
The bundles $U_i$ over $S_i$  are trivial for $i=0,1$. 
\end{lem}

\begin{proof}
For  $i=0$ case, the bundle $U_0$ is trivial as any bundle over a point is trivial. 

For   case $i=1$ case,   the bundle $U_1$ is a rank-$(d-1)$  subbundle of a trivial bundle 
$S_1 \times \A^d \to S_1$ whose fibre over ${\bf w}$ consists of the quotient space 
$\A^{d} / \langle {\bf w} \rangle$.  On the other hand, the  subbundle $N_1$ which is ``normal" to 
$U_1$, that is having the affine line $\A^1 =  \langle {\bf w} \rangle $ spanned by ${\bf w}$  
as the fibre over ${\bf w}$ is trivial:
The map 
\begin{eqnarray*}
 S_1  \times \K   \to    N_1: \ \ \ 
(  {\bf w}   ,   c )  \mapsto   (  {\bf w} ,  c \cdot {\bf w} ) 
\end{eqnarray*}
is an isomorphism and provides the trivialization that is needed. This implies that 
$U_1$ is also a trivial bundle over $S_1$. 
\end{proof}

\subsubsection{The class of the simple Feynman quadric   in the Grothendieck ring}

\begin{lem}
The class $[ {\cQ}]$  of the quadric in the Grothendieck ring $K_0(Var)$ is given by
$$
[ {\cQ}] = [  U_0] + [  U_1]  = \bL^{2d-1}  + \bL^d -  \bL^{d-1} 
\in \Z [\A^1] \subset K_0(Var)
$$  
where $\bL$ denotes the Leftschetz motive $[\A^1]$. 
\end{lem}

\begin{proof}
The scissor congruence allows us write the class $[ {\cQ}] $ as the sum 
$ [  U_0] + [  U_1] $, see \eqref{eqn_cell_fibration}. 
The rest follows from the fact that, for  a locally trivial fibration $E \to B$ with fibre 
$F$, the classes in the Grothendieck ring of varieties  $K_0(Var)$ satisfies
$[E] = [B] \cdot [F]:$
\begin{eqnarray*}
 [  U_0] + [  U_1]   
&=&   [\A^0] \cdot  [\A^d]  +  ( [\A^d] -1 )  \cdot  [\A^{d-1}]   \\
&= &  \bL^{2d-1}  + \bL^d -  \bL^{d-1}.
\end{eqnarray*}
\end{proof}

\begin{rem} \label{rem_affine}
The  same line of arguments shows that the class $[\widehat {\cQ}]$ of the  quadric 
in the Grothendieck ring $K_0(Var)$ is   $  \bL^d + [\P^{d}] \cdot [\P^{d-1}]  $: Consider 
the restriction $\widehat\pi: \widehat \cQ \subset \P^d \times \P^d \to  \P^d$ of the projection.    
Then, the  fibres are
 \begin{eqnarray} \label{eqn_affine_fibre}
\widehat {\pi}^{-1}( [ {\bf w} ]  ) \cong 
\left\{
	\begin{array}{lcl}	
		\P^{d} = \{(  [ z_0 : \dots : z_d ] \}        
					&	\text{if} & [{\bf w}] = [1: 0 : \dots : 0],  \\
		\P^{d-1} = \{ [ z_0 : \dots : z_d ]  \mid  z_1w_1 + \cdots + z_d w_d  =0 \}
					&	\text{if} &  [{\bf w}] \ne  [1: 0 : \dots : 0]. 		
	\end{array}	
 \right.  
\end{eqnarray}
Therefore, we have a decomposition according to the types of the fibres  in \eqref{eqn_affine_fibre}:
\begin{equation*}
  \widehat S_{-1} = \emptyset, \ \ \  \widehat  S_0 = \{  [1: 0 : \dots : 0] \} \subset \P^d ,  \ \ \  
  \widehat  S_1  = \P^d \setminus  \widehat  S_0,
\end{equation*}
and $\widehat U_i = \widehat {\pi}^{-1}( \widehat  S_i), i=0,1$. 
 The scissor congruence and the classes of  fibrations in $K_0(Var)$  provides that 
\begin{eqnarray*}
 [\widehat U_0] + [\widehat U_1]  
&=&   [\A^0] \cdot  [\P^d]  +  ( [ \P^{d}]  -1 )  \cdot  [\P^{d-1}]   \\
&= &  ( \bL^d +  [ \P^{d-1}] )  +   (\bL + \cdots + \bL^d)  \cdot    [\P^{d-1}]  \\
&= &  \bL^d +  (1 + \bL + \cdots + \bL^d)  \cdot    [\P^{d-1}]   = \bL^d +  [\P^{d}] \cdot [\P^{d-1}].
\end{eqnarray*}
 \end{rem}

\subsubsection{The motive of the simple Feynman quadric   in Voevodsky's category}

\begin{lem}
\label{lem_mot_fib}
The motives $\m^c(  U_i), i=0,1$    associated to the fibrations  \eqref{eqn_cell_fibration} 
are given by are
\begin{eqnarray*}  
\m^c(  U_i)  
& =  & 
	\left\{ 
		\begin{array}{lc}
		 \m^c(S_0 \times \A^d) = \m^c( \A^0)(d) [2d]   \ \ \ \& \ \ \ i = 0 \\
 		  \m^c(S_1 \times \A^{d-1})  =  \m (\A^{d} \setminus \{0\})(p-1)[2p-2] \ \ \  \& \ \ \ i = 1. \
 		\end{array}
	\right. \\
\end{eqnarray*}
\end{lem}

\begin{proof}
From  Corollary 4.1.8 of \cite{Vo} we know that taking the product with an affine space $\A^k$ is an isomorphism 
at the level of the corresponding motives: the motive $\m^c(X \times \A^k) = \m^c(X)(k) [2k]$ is obtained from  
$\m^c(X)$ by Tate twists and shifts. 
The result then follows by applying this identity to the fibrations \eqref{eqn_cell_fibration}.
\end{proof}


\begin{cor} \label{cor_mot}
The motives $\m^c(  U_i), i=0,1,2$     are mixed Tate motives. 
\end{cor}

\begin{proof}
 This is an immediate consequence of Lemma \ref{lem_mot_fib}: The motives  $\m^c(  U_i$)  of
the fibrations  depend only on the motives   $\m^c(  S_i$) of the base and   the motives   
$\m^c(F_i$) of the the fibres  through products, 
Tate twists, sums, and shifts. All these operations preserve the subcategory of mixed Tate 
motives. The result follows from the fact that the motives of $  S_i$ and  
$  F_i$  are mixed Tate. 
\end{proof}

\begin{prop} 
\label{prop_motive1}
The Voevodsky motive $\m^c( \cQ)$     associated to the quadric $ \cQ$ is 
a mixed Tate motive. 
\end{prop}

\begin{proof}  
The bundle $  U_0  $ is a closed subscheme in $ \cQ$. We can
use the canonical distinguished triangle (Prop 4.1.5 in \cite{Vo}) that is
\begin{eqnarray*}  
\begin{diagram}
\m^c(  U_0  )   & \rTo &  \m^c( \cQ) & \rTo & 
		\m^c(  \cQ \setminus U_0  )
		& \rTo &   \m^c(  U_0  )[1].  
\end{diagram}
\end{eqnarray*}
The motives  $\m^c(  U_0 )$ and $\m^c(  \cQ \setminus U_0  ) = \m^c(U_1)$ are mixed Tate
due to  Corollary \ref{cor_mot} which simply implies that the motive  $\m^c( \cQ) $ of the
simple Feynman quadric is also mixed Tate. 
\end{proof}

\begin{rem}  
The  same line of arguments shows that the motive $\m(\widehat {\cQ})$ of the Zariski 
closure $\widehat {\cQ}$ of  quadric inside $\P^d \times \P^d$ 
 also defines  a mixed Tate motive. The main difference here is that 
the fibres in \eqref{eqn_affine_fibre} are no longer affine spaces and
the product formula is usable.

However, the results in \cite{Habibi} provides usable setup such cases: The motives of the 
proper smooth locally trivial fibrations $E \to B$ with fibres $F$ that admit cellular 
decomposition and satisfy  Poincar\'e  duality can be given as 
$$
\m(E) =  \bigoplus_{p \geq 0} CH_p(F) \oplus \m (B)(p)[2p]
$$
due to Thm 2.10 in \cite{Habibi}.
\end{rem}

\subsection{Singularity of the simple  Feynman quadric, its resolution and its motive}
\label{sec_singular}

The   vector 
$$ 
(\frac{\partial {  q}}{\partial {\bf z}},\frac{\partial {  q}}{\partial {\bf w}}) := 
(\frac{ \partial q }{\partial z_1 }, \dots , \frac{ \partial q }{\partial z_d} ; \frac{ \partial q }{\partial w_1} ,  \dots ,  \frac{ \partial q }{\partial w_d } )
 = (w_1, \dots, w_d; z_1, \dots ,w _d ) 
\in T_{({\bf z} ; {\bf w})} (\A^d \times \A^d)
$$
that is normal to $\cQ$ at $({\bf z} ; {\bf w}) \in \cQ$ does not vanish unless $({\bf z} ; {\bf w}) = (0;0)$, i.e.,  $\cQ$ has no singular points other
than $(0;0)$.  It is easy check that the same is true for $\widehat\cQ \setminus (0;0)$. 
We denote the smooth part $\cQ \setminus (0;0)$ 
(resp. $\widehat\cQ \setminus (0;0)$) by $\cQ_{sm}$ (resp. by $\widehat\cQ_{sm}$).

The singularity of $\cQ$ is a quite simple type: 
Our affine quadric $ {\cQ} $ admits the following multiplicative group $\mathbb{G}_m$ action;
$$
\forall \lambda \in \mathbb{G}_m, \ \ \ (z_1,\dots,z_d ; w_1,\dots,w_d) \mapsto  
(\lambda z_1,\dots, \lambda z_d  ; \lambda w_1,\dots,\lambda w_d,)
$$
that provides us the   fibration
\begin{equation} \label{eqn_projectivization}
\begin{diagram}
	\mathbb{A}^1 \setminus \{0\}	& \rInto	& {\cQ}_{sm} \\
	\dTo						&		& \dTo \\
	\text{pt}					& \rInto	& {\cQ}_{sm} / \mathbb{G}_m.
\end{diagram}
\end{equation}
In other words, $\cQ$ can be thought as  the cone of the smooth  quadric 
$$
{\cQ}_{sm} / \mathbb{G}_m = \{ [z_1,\dots,z_d ; w_1,\dots,w_d] \mid 
z_1 w_1 + \cdots + z_d w_d  = 0\} \subset \P^{2d-1}.
$$ 
Moreover, it also  hints us that the motive of ${\cQ}_{sm} / \mathbb{G}_m$ is mixed Tate.

\begin{lem}
$[ {\cQ}_{sm} / \mathbb{G}_m ]  = (\bL^{2d-1} +\bL^d - \bL^{d-1} -1) \cdot   \left(\sum_{r = 0}^\infty - \bL^r \right)$.
\end{lem}

\begin{proof}
As we know that the class of the fibres $\bL -1$ and the class $[\cQ_{sm}] = [\cQ] -1$  of the total space (see,
Remark \ref{rem_affine}),  the class 
$[ {\cQ}_{sm} / \mathbb{G}_m ]  $
of the quotient can be given by $ [\cQ_{sm}] \cdot  (\bL -1)^{-1}$  as stated above. 
\end{proof}

\subsubsection{Resolution of the singularity}
 The morphism
\begin{eqnarray*}
\sigma: \A^{2d} & \to &  \A^{2d} \times \P^{2d-1} \\
(z_1,\dots,z_D,w_1,\dots,w_D) & \mapsto & ((z_1,\dots,z_D,w_1,\dots,w_D), [z_1 : \dots: z_D : w_1 : \dots : w_D])
\end{eqnarray*}
 is defined outside the origin of $\A^{2d}$. The  Zariski-closure $ \text{Bl}_{(0;0)} {\cQ}$  of  
$\sigma( \cQ_{sm}) $ inside $\A^{2d} \times \P^{2d-1} $ is a nonsingular subvariety. It
is actually blow-up  of $ {\cQ}$  resolving the singularity at $(0;0)$.
We denote the exceptional divisor $( \text{Bl}_{(0;0)} {\cQ} ) \setminus \cQ_{sm})$ by $\cD$.

\begin{prop} \label{prop_blow}
The  motive $ \text{Bl}_{(0;0)} {\cQ}$  of  the blown-up affine quadric  $\cQ$  is a mixed Tate motive. 
\end{prop}

\begin{proof}
The strict transform is an isomorphims away from the singular point
of $ {\cQ} $. We show that $\m^c(\cQ_{sm})$ is mixed Tate by simply using the canonical
distinguished triangle for $(0;0) \hookrightarrow {\cQ}$.

We only need to calculate the motive of the exceptional divisor $\cD$ and patch them together using the
  distingushed  triangle for $\cD \hookrightarrow \text{Bl}_{(0;0)} {\cQ}$. 
Consider the blow-up  $ \text{Bl}_{(0;0)} \A^{2d}$ at $(0;0)$  as a subspace of $\A^{2d} \times \P^{2d -1}$
with the projection $\A^{2d} \times \P^{2d -1} \to \A^{2d}$. 
If we pick a homogeneous chart  $[u_1 : \dots: u_d : v_1 : \dots : v_{d}]$ 
in the central fibre  $\{0\} \times  \P^{2d -1}$ of  the projection,
then, we observe that 
$$
\cD = \text{Bl}_{(0;0)} {\cQ} \cap ( \{0\} \times \P^{2d -1} ) = \{ [u_1 : \dots: u_d : v_1 : \dots : v_{d}] \mid u_1 v_1 + \cdots   +  u_d v_d=0\}. 
$$

We can give a motivic  decomposition of $\cD$    
as in \S \ref{sec_projection}. The  same line of arguments  in the proof of the Proposition \ref{prop_motive1} 
applies to this stratification and   implies that  the exceptional divisor
is of type mixed Tate.  
\end{proof}


\begin{cor} 
The  motive $\m^c(\widetilde{\cQ})$  of  the blow-up of the quadric $\widehat\cQ$ at its singular point  is also a mixed Tate motive. 
\end{cor}

\begin{proof}
This statement follows from the facts that the motives $\m^c(\text{Bl}_{(0;0)} {\cQ})$
and $\m^c(\widehat \cQ \setminus {\cQ})$ of the strata o   $\widetilde{\cQ}$ 
are mixed Tate. 
\end{proof}

\begin{rem}  \label{rem_homotopy}
Proposition \ref{prop_blow} can be proved alternatively as follows.
The fibration  
\eqref{eqn_projectivization} is  modified by adding the points at infinity to each fibre 
$\A^1 \setminus \{0\}$, and that provides an locally trivial $\A^1$-fibration:
\begin{equation} 
\begin{diagram}
	\A^1 \cong \ 	& \mathbb{P}^1 \setminus \{0\}	& \rInto	& \widehat {\cQ}_{sm}	&  \ =  \widehat {\cQ} \setminus \{0\} \\
				& \dTo					&		& \dTo				& \\ 
				& \text{pt}					& \rInto	& {\cQ}_{sm} / \mathbb{G}_m. &
\end{diagram}
\end{equation} 
We can conclude that  the motive  $\m(  {\cQ}_{sm} /  \mathbb{G}_m)$  of the exceptional 
divisor  is mixed Tate as  the motives of the locally   trivial $\A^1$-fibrations are 
the same as their bases due to $\A^1$-homotopy invariance. 
\end{rem}

  


\section{The position space Feynman quadrics}  
\label{sec_quadric}

This section is the technical heart of our paper. We associate  a  Feynman quadric
to each  Feynman graphs and study them, in particular, in the case of the complete graphs. 
The reduction to the complete graphs 
can justified by the fact the periods of Feynman quadrics can be formulated
as the periods of  the Feynman quadric of  complete graph with same number of vertices 
(see \cite{CeyMar3} for a similar treatment for configuration space setup). We
simply imitate the projections in \S \ref{sec_projection} and introduce the complement
of the Feynman quadric as a configuration space of certain hyperplane arrangements. Then,
we give a stratification of these configuration spaces in terms of the degeneration types
of these hyperplane arrangements. 


\subsection{Feynman quadric associated to Feynman graphs} 
\label{sec_feynman_quadric}

The {\it  Feynman quadric} $Z_\Gamma$ associated to a given Feynman graph $\Gamma$ is the 
 quadric 
$$
Z_\Gamma :=  \bigcup_{ e \in Edg(\Gamma) } \cH^{e}  \subset (\A^d \times \A^d)^{Ver(\Gamma)}
$$
whose irreducible components are 
$$
\cH^{e} := \{ q^{ij} = q( {\bf z}^i - {\bf z}^j ,  {\bf w}^i - {\bf w}^j  ) = 0 \mid (ij) = \partial_\Gamma (e) \}.
$$

\begin{lem} \label{lem_quadric}
Let $e \in Edg(\Gamma)$ and $\partial_\Gamma (e) = (ij)$. The quadric 
$\cH^{e} \subset (\A^d \times \A^d)^{Ver(\Gamma)}$ is isomorphic to 
$\cQ \times (\A^d \times \A^d)^{Ver(\Gamma) \setminus \{j\}}$. 
\end{lem}

\begin{proof} 
Consider  the composition of the morphisms:
\begin{eqnarray*}  
\begin{diagram}
(\A^d \times \A^d)^{Ver(\Gamma)} 	&	 \rTo^{p_{ij}}	& 	(\A^d \times \A^d)^{\{i,j\}} 	& \rTo^{tr_i}	&	\A^d \times \A^d  \\
({\bf z}^m ;  {\bf w}^m \mid m \in Ver(\Gamma))  &  \rMapsto &  ( ({\bf z}^i ; {\bf w}^i) , ({\bf z}^j ; {\bf w}^j )  )  &  \rMapsto &  
( {\bf z}^i - {\bf z}^j ;  {\bf w}^i - {\bf w}^j   ).
\end{diagram}
\end{eqnarray*}
The  morphism $p_{ij}$ simply forgets all factors but $ ({\bf z}^i ; {\bf w}^i)$ and  $({\bf z}^j ; {\bf w}^j ) $, hence a trivial 
fibration with fibres $(\A^d \times \A^d)^{Ver(\Gamma) \setminus  \{i,j\}}$. On the other hand the morphism $tr_i$uses 
the translations
to fix  $ ({\bf z}^i ; {\bf w}^i)$ at $(0;0)$, therefore it is also a trivial fibration with fibres $\A^d \times \A^d$ parameterizing  
$ ({\bf z}^i ; {\bf w}^i)$. 
The quadric $ {\cH}^{e}$ is simply the preimage $ (tr_i \circ p_{ij})^{-1}  (\cQ)$ of the simple Feynman quadric $\cQ$. 
\end{proof}

\begin{cor} \label{cor_diagonal}
The singular locus of  the quadric $\cH^{e} $ is the diagonal
$$
\Delta^{e} = \{ ({\bf z}^m,  {\bf w}^m  \mid m \in Ver(\Gamma)) \in (\A^d \times \A^d)^{Ver(\Gamma)}  \mid {\bf z}^i = {\bf z}^j \ \& \  {\bf w}^i = {\bf w}^j   \}.
$$
\end{cor}

\begin{proof}
The singularities of  $\cH^{e} $ 
are determined as stated by using the singularities of $\cQ$ 
in \S \ref{sec_singular} and Lemma \ref{lem_quadric}.
\end{proof}

\subsection{A projection for the Feynman quadrics } 
In this paragraph, we  establish an explicit connection between the (complements of the) 
Feynman quadrics and  a  configuration spaces of  hyperplane arrangements in ``almost general
position".


\subsubsection{Projection and  its fibres} 
\label{sec_projection_extended}

This paragraph is a straightforward  generalization of \S  \ref{sec_projection} which essentially 
examines   the case of one-edge graph.

For any $\Gamma$ with $|Ver(\Gamma)| = n$, the Feynman quadric $Z_\Gamma$ is  contained (set theoretically)
in $Z_{\kappa_n}$ where  $\kappa_n$ is the  complete graph with $n$-vertices. From now on,  we will 
consider only  the complete graphs. However, the main strategy below is valid for any Feynman 
graphs $\Gamma$. We will remark the ramifications in \S \ref{sec_corollaries} after completing case of complete graphs.

\begin{figure}
\begin{center}
\includegraphics[scale=1.1]{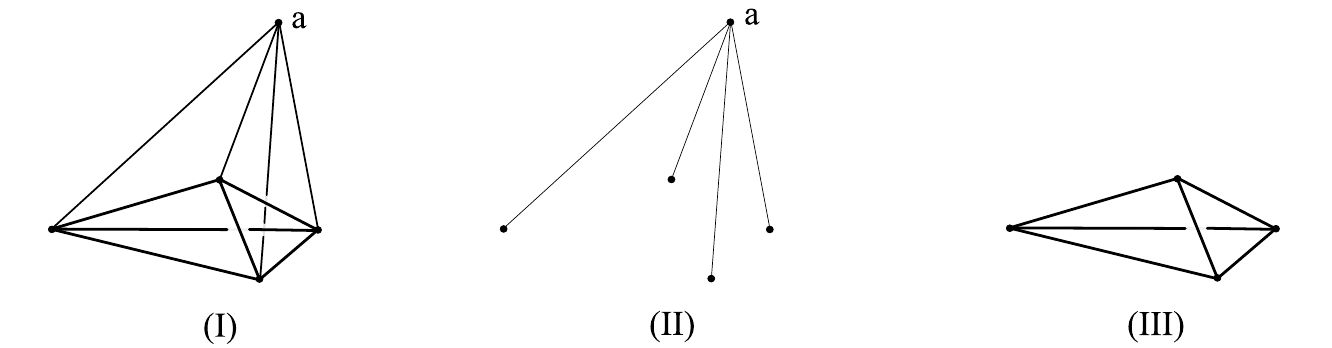}
\caption{(I) $\kappa^+$ complete graph with 5 vertices; (II) Star of $a$; (III) $\kappa$, obtained from $\kappa^+$ be removing $St(a)$. 
\label{fig_graph}}
\end{center}
\end{figure}

Let $\kappa^+, \kappa$ be a pair of complete  graphs such that $\kappa$ is obtained from
$\kappa ^+$ by removing one of its  vertices  $a \in Ver(\kappa ^+)$, and the all edges 
adjacent to the vertex $a$, i.e.,  $Ver(\kappa ^+) = Ver(\kappa) \cup \{a\}$, $Edg(\kappa ^+) = Edg(\kappa) 
\cup St(a)$ (which is called {\it the star of the vertex $a$}) where $St(a)$ denotes the 
set $\partial^{-1}_{\kappa ^+} (a)$ of edges adjacent to the vertex $a$,
and the boundary map $\partial_{\kappa}$ is the restriction of $\partial_{\kappa ^+}$ to
$Edg(\kappa)$.   For an example, see Figure \ref{fig_graph} which  illustrates the simple case where $\kappa^+$ is  
complete graph with 5-vertices, the star $St(a)$ $a$ and $\kappa$.  

Consider the projection 
\begin{equation} 
	\begin{array}{rcl} 	\label{eqn_full_forgetful-1}
	\pi_{\kappa ^+}: (\A^d \times \A^d)^{Ver(\kappa ^+)}		& \to		& (\A^d \times \A^d)^{Ver(\kappa)} \times \A^d   \\ 
	(  {\bf z}^m  ,  {\bf w^m}  \mid m \in Ver(\kappa ^+) )		&	\mapsto	& 
										{\bf b} := (  {\bf z}^m  ,  {\bf w^m}  \mid m \in Ver(\kappa) ) \times  {\bf w}^a
	\end{array}
\end{equation}
whose fibers are  $\P^d$. 

In order to examine the Feynman quadric  
$ Z_{\kappa ^+} = \bigcup_{e \in Edg(\kappa ^+)}\cH^{e} \subset (\A^d \times \A^d)^{Ver(\kappa ^+) }$ for
a bigger graph $\kappa^+$, 
we study the fibres of the restriction of the forgetful morphism \eqref{eqn_full_forgetful-1}:
\begin{equation} 
	\begin{array}{rcl}	\label{eqn_full_forgetful}
	 \pi_{\kappa ^+}:    Z_{\kappa ^+} 	& \to		& (\A^d \times \A^d)^{Ver(\kappa)} \times \A^d .
	\end{array}
\end{equation}

\begin{lem} \label{lem_full_fibre}
The  fibre $\pi^{-1}_{\kappa ^+} ( {\bf b} )$   of    \eqref{eqn_full_forgetful}
over a point  ${\bf b} $   is
\begin{enumerate}
	
	\item the affine space $\A^d$, if ${\bf b}  \in (Z_{\kappa}  \times \A^d)$, 
	
	\item the affine space $\A^d$, if $ {\bf b}  \in   \bigcup_{e \in St(a)  } \{   {\bf w}^a  =  {\bf w}^i  \mid \partial_{\kappa ^+}(e)= (ai) \} 
			\subset ((\A^d \times \A^d)^{Ver(\kappa)} \times \A^d)$, 
	
	\item a    hyperplane arrangement  
	$$
	\cA_{St(a)} := (  P_{e} \in \A^d \mid   e \in St(a) )  	
	$$ 
	 where 
	$$
	P_{e} :=\{ {\bf z}^a   \in \A^d  \mid q^{ai} = q({\bf z}^a - {\bf z}^i, {\bf w}^a - {\bf w}^i) = 0 \ \text{and} \ (ai) = \partial_{\kappa ^+} (e) \}, 
	$$
	when ${\bf b} \not\in (  Z_{\kappa}  \times \A^d) \cup \bigcup_{e \in St(a)) } \{   {\bf w}^a  =  {\bf w}^i  \}$.
	
	\noindent
	Moreover, for any $I \subset St(a)$ with $|I| \leq d$, 
	 the subarrangements $\cA_I = (P_{e}  \mid e \in I)$ satisfy
	\begin{equation*}
	P_{I}  :=  \bigcap_{{e \in I }   } P_{e} = \A^{d - |I|}.
	\end{equation*}	 
	 In other words, they are in general position.

%
%
%
\end{enumerate}
\end{lem}

\begin{proof}
The statement directly follows from the ideals  defining of the irreducible
quadrics $\cH^{e}$:

(1)  If ${\bf b}  \in (Z_ \kappa  \times \A^d) \subset ((\A^d \times \A^d)^{Ver(\kappa)} \times \A^d)$, then the equation 
$$
\prod_{ {e \in Edg(\kappa ^+)} \atop {\partial_{\kappa ^+} (e) = (ij)}} q^{ij} =  
\left(\prod_{ {e \in Edg(\kappa)} \atop {\partial_{\kappa} (e) = (ij)} } q^{ij} \right) \cdot 
\left( \prod_{ {e \in St(a) } \atop {\partial_{\kappa ^+} (e) = (ai )}} q^{ai} \right) =
0 \cdot \left( \prod_{ {e \in St(a) } \atop {\partial_{\kappa ^+} (e) = (ai )}} q^{ai} \right)  =0   
$$ 
is satisfied for any ${\bf z}^a$
in the fibre  $\A^d =  \pi^{-1}_{\kappa ^+}  ({\bf b})$ of   \eqref{eqn_full_forgetful-1}. 
Therefore, all  ${\bf z}^a \in \A^d $ must be  in $ Z_{\kappa ^+}$.

(2)  If ${\bf b} \in \{   {\bf w}^a  =  {\bf w}^i \}$ for an edge  $e \in St(a)$ with $\partial_{\kappa ^+}(e)= (ai)$, then the equation 
$q^{ai} =  ({\bf z}^a - {\bf z}^i) \cdot ({\bf w}^a - {\bf w}^i) =  ({\bf z}^a - {\bf z}^i) \cdot 0 =0$ is satisfied for any ${\bf z}^a$
in the fibre  $\A^d =   \pi^{-1}_{\kappa ^+}  ({\bf b}) $  of   \eqref{eqn_full_forgetful-1}.
Therefore, all  ${\bf z}^a  \in \A^d $ must be  in   $   \cH^{e} \subset     Z_{\kappa ^+}$ in
such a case.

(3) In all other cases, the intersections  of  $ \cH^{e}$  for $e \in St(a)$ with the fibres 
$\A^d =  \pi^{-1}_{\kappa ^+}  ({\bf b})$   of the forgetful morphism \eqref{eqn_full_forgetful-1}
are simply defined by the equations $q^{ai} = ({\bf z}^a - {\bf z}^i) \cdot  ({\bf w}^a - {\bf w}^i) =0$ where 
$\partial_{\kappa ^+}(e) = (ai)$. Therefore, the fibre
$(\bigcup_{e \in St(a)}  \cH^{e}) \cap  \pi^{-1}_{\kappa ^+}  ({\bf b})$ of \eqref{eqn_full_forgetful}  
over ${\bf b}$ is the 
hyperplane  arrangement as stated. Note that, $q^{ai} = 0$ cannot be  
the hyperplane infinity, since that hyperplane is given by the equation $z^a_0=0$.
 
 We only need to show that  the  subarrangements $  \cA_I = (  P_{e}  \mid e \in I)$ 
 is in general position for all $I \subset St(a)$ with $|I| \leq d$ .  

As we consider the compliment of cases examined 
 in (1) and (2), we can simply set  
$$
P_e = \{ {\bf z}^a  \mid ({\bf z}^a - {\bf z}^i) \perp  ({\bf w}^a - {\bf w}^i) \} . 
$$ 
Note that, if a subarrangement $\cA_I$  of $ \cA_{St(a)} $   is in general position,
then the  vectors ${\bf v}^{i} = {\bf w}^{a} - {\bf w}^{i} $  that are normal to the  affine   hyperplanes 
$P_{e} =  \{ q^{ai} =  ({\bf z}^a - {\bf z}^i) \cdot  ({\bf w}^a - {\bf w}^i) =0 \mid \partial_{\Gamma^+}(e) = (ai) \} $ for $e \in   I$ must 
span a $|I|$-dimensional vector space when $|I| \leq d$, or simply $d$-dimensional 
vector space when $|I| \geq d$ . We prove the statement by induction
on the cardinality of $I \subset St(a)$:

First step, $|I|=2$ case: Let  $I= \{e_i, e_j\}$ and $\partial_{\kappa ^+}(e_*) = (a*)$. 
If the statement does not hold  for the pair  $   P_{e_i},   P_{e_j}$,   
then the affine part of these hyperplanes  must be parallel. Hence their normal vectors 
satisfy 
\begin{equation} \label{eqn_parallel} 
({\bf w}^a - {\bf w}^i) = \lambda  ({\bf w}^a - {\bf w}^j)
\end{equation}
 for a nonzero $\lambda$. In this case, we have
\begin{equation} \label{eqn_ohlala} 
	\begin{array}{ccc}
		\begin{array}{c}
		({\bf z}^a - {\bf z}^i) \perp  ({\bf w}^a - {\bf w}^i) \ \&  \\
		({\bf z}^a - {\bf z}^j) \perp  ({\bf w}^a - {\bf w}^j)
		\end{array}
	& 
		\begin{array}{c}
		\Longrightarrow  \\
		\text{due to} \eqref{eqn_parallel}
		\end{array}
	&
		\begin{array}{c}
		({\bf z}^a - {\bf z}^i) \perp  ({\bf w}^a - {\bf w}^i)  \ \& \\
		({\bf z}^a - {\bf z}^j) \perp  \lambda ({\bf w}^a - {\bf w}^i)
		\end{array}
	\end{array} \
\end{equation}
which implies that $({\bf z}^a - {\bf z}^j)  - ({\bf z}^a - {\bf z}^i) =  ({\bf z}^i - {\bf z}^j) \perp  ({\bf w}^a - {\bf w}^i)$.
After interchanging indices $i$ and $j$, the same argument   implies 
$({\bf z}^i - {\bf z}^j) \perp  ({\bf w}^a - {\bf w}^j)$. We conclude that 
$$
	\begin{array}{ccc}
		\begin{array}{c}
		({\bf z}^i - {\bf z}^j) \perp  ({\bf w}^a - {\bf w}^i) \\
		({\bf z}^i - {\bf z}^j) \perp  ({\bf w}^a - {\bf w}^j)
		\end{array}
	&
		\begin{array}{c}
		\Longrightarrow  
		\end{array}
	&
	({\bf z}^i - {\bf z}^j) \perp ({\bf w}^i - {\bf w}^j) = ({\bf w}^a - {\bf w}^j) -  ({\bf w}^a - {\bf w}^i).
	\end{array} \
$$
In other words, the normal vectors $({\bf w}^a - {\bf w}^i)$ and $({\bf w}^a - {\bf w}^i)$ can be parallel
only when ${\bf b} \in  \cH^{e} \times \A^d \subset   Z_ \kappa \times \A^d$.
That  contradicts with our initial assumption on ${\bf b}$.


Next, we assume that the statement  holds  for $I = \{e_{i_1},\dots,e_{i_k} \} $ with $|I| < d$ but not for 
$J = I \cup \{e_{i_{k+1}} \} \subset St(a)$. 
Then, the sets of normal vectors ${\bf v}^{i_1} , \dots , {\bf v}^{i_{k-1}} , {\bf v}^{i_{k}}$ and
 ${\bf v}^{i_1} , \dots , {\bf v}^{i_{k-1}}, {\bf v}^{i_{k+1}}$ must span the same vector space. This 
 can be satisfied if and only if
 $$
( {\bf v}^{i_1}  \wedge \dots \wedge  {\bf v}^{i_{k-1}} ) \wedge {\bf v}^{i_{k}} =
\lambda   \cdot ( {\bf v}^{i_1}  \wedge \dots \wedge  {\bf v}^{i_{k-1}}  ) \wedge {\bf v}^{i_{k+1}}  
 $$
 for a nonzero $\lambda$. This equality implies  that $ {\bf v}^{i_{k}}  =  \lambda  \cdot  {\bf v}^{i_{k+1}} $,
 which again contradicts with the assumption that ${\bf b} \not \in   Z_ \kappa \times \A^d$.
 Therefore, the arrangement $\cA_{J} = (  P_{e}  \mid e \in J)$  with  $|J| \leq d$ must be also be 
 in general position and,  can be characterized in terms of intersections as stated. 
\end{proof}

\subsection{Configuration space of hyperplane arrangements in ``almost general position"} 
\label{sec_conf_space}

The case (3)  in Lemma \ref{lem_full_fibre} suggests us that  the space
$$
\cM^{St(a)}  := 
((\A^d \times \A^d)^{Ver(\kappa)} \times \A^d) \setminus ((  Z_ \kappa \times \A^d) \cup 
\bigcup_{e \in St(a)} \{  {\bf w}^i  =  {\bf w}^a  \mid (ai) = \partial_{\kappa} (e) \})
$$
can be thought as a configuration space of  the  hyperplane arrangements satisfying certain 
conditions, i.e., 
it parameterizes the arrangements  $  \cA_{St(a)} = (  P_{e}  \mid e \in St(a)) $  
of hyperplanes in $\A^d$  labeled by the index set $St(a)$ such that the intersections
\begin{equation} \label{eqn_generic}
  P_I := \bigcap_{e \in I}   P_{e} =
		\A^{d - |I|}		\ \ \ \ \text{for each} \ \ I \subset St(a) \ \ \text{with} \ \ |I| \leq d.
\end{equation}
We call these hyperplane arrangements are in {\it almost general position} as they can be put
into general position by using the parallel translations of (at most $|St(a)|-d$) hyperplanes. 
\begin{figure}
\begin{center}
\includegraphics[scale=1.1]{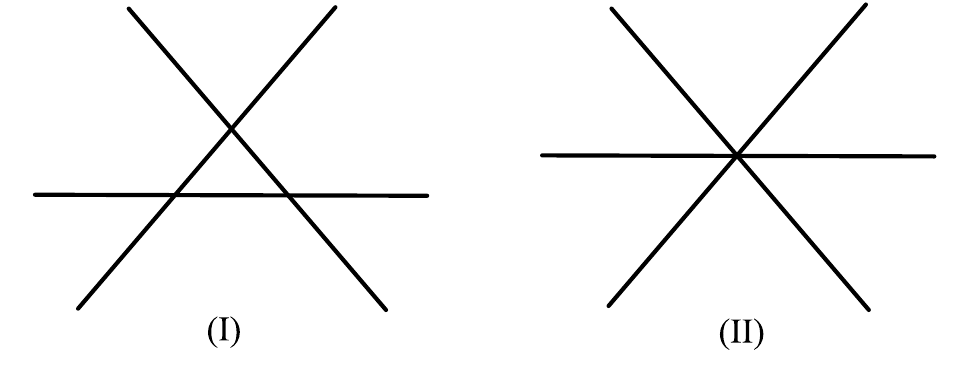}
\caption{Arrangements of 3 hyperplanes in $\A^2$: (I) In general positions, (II) In almost general 
(but not in general) position. \label{fig_triple}}
\end{center}
\end{figure}


\subsubsection{Degeneration types of hyperplane arrangements in  almost general position} 
\label{sec_degeneration}

Matroids and intersection posets  provide  general setups for encoding the 
degeneration types of the hyperplane arrangements, see, for instance \cite{Alexeev,Stan}. As we 
only need to consider the degenerations of hyperplane arrangements  of very particular type,  we
 provide a simpler setup below.

Here, we note that a part  of the conditions imposed by \eqref {eqn_generic} on the 
{\it intersection poset} $L(  \cA_{St(a)})$ is implicit. 

\begin{lem} \label{lem_degeneration_types}
Let  $  \cA_{St(a)}$ (resp. $  \cA_{St(a)}'$) be a hyperplane arrangement in (resp. almost) 
general position. Then,  the map
\begin{equation}
\label{eqn_bijection}
L(  \cA_{St(a)}) \to L(  \cA'_{St(a)}):    P_I \mapsto   P'_I
\end{equation}
is injective,  and the complement of its image   $\{   P'_J \mid J \in \Theta' \}$ is indexed by the set
\begin{equation}
\label{eqn_degeneration_type}
\Theta' := \{ J \subset St(a)  \mid  \ |J| > d \ \& \  {P}'_J  =  \{ point \}  \} \subset L(  \cA'_{St(a)})
\end{equation}
of deepest intersections.
\end{lem}

We call the possible index sets \eqref{eqn_degeneration_type}   as the  {\it  degeneration types} of the 
hyperplane arrangements   in almost general position.
If $\Theta' =   \emptyset $, then  $L(  \cA_{St(a)}) = L(  \cA'_{St(a)})$ and the  
hyperplane arrangement $  \cA'_{St(a)}$ is also in  general position. 
If $\Theta'  \ne   \emptyset $, we call $  \cA'_{St(a)}$ a   {\it degenerate} hyperplane 
arrangement.

\begin{proof} (of Lemma \ref{lem_degeneration_types})
For the subsets $J$ with $| J | \leq d$,  the condition \eqref {eqn_generic} is explicit and implies
that the intersections $  P_J,    P'_J$  are elements respectively in $  \cA_{St(a)}$ 
and  $  \cA'_{St(a)}$. Hence,  the map \eqref{eqn_bijection}
is bijection for the index set $\{J \subset St(a) \mid | J |  \leq d \}$. 

Since $ \cA_{St(a)}$ is in general position,   the deeper intersections $  P_J = \emptyset$ for $J \in L(  \cA_{St(a)})$  with
$| J | > d$, and we conclude that the difference between the  intersection posets of these hyperplane 
arrangements can be given  by the following index set  
$$
\{ J \subset St(a)  \mid   |J| > d \ \& \  {P}'_J \not= \emptyset \}
$$

However,  such  intersections $  P'_{J} \in L(  \cA'_{St(a)})$
are not arbitrary and   being in almost general position \eqref{eqn_generic}  implies certain implicit  conditions 
on $  P'_{J}$.  For instance, while the triple intersections in Figure \ref{fig_triple}
in $\A^2$ provides arrangements in almost general position, a similar triple intersection in $\A^3$ does not. 
More precisesly, 
there may exists hyperplanes  arrangement  in almost general position in 
$\A^d$  with more then $d$  hyperplanes intersect at the same point.  
However, note also that, an intersection along a higher dimensional subspace for  $|J| > d$ is prohibited  by the   
condition \eqref{eqn_generic}: If there exists $  P'_J = \A^k, k >0$ for a subset $|J| > d$,
then there must be a subset $K \subset J$ with $d \geq |K| > d -k$, such that, the corresponding
subarrangement $(  P_{e}  \mid e \in K) $ of hyperplanes violate \eqref{eqn_generic}, i.e.,
$ {P}'_J$ can only be a single point as stated.
\end{proof}
 

The following stratification of the configuration space $\cM^{St(a)}$ is a tautology:  

\begin{prop} 
\label{prop_stratification}

\begin{enumerate}
\item  For any given degeneration type $\Theta$,  there is a quasi-projective subvariety
$S_\Theta \subset \cM^{St(a)}$ parameterizing the hyperplane arrangements in almost general position 
with the    degeneration type   $\Theta $. 

\item The configuration space $\cM^{St(a)}$ is stratified by these pairwise disjoint subvarities 
$S_\Theta$. 
\end{enumerate}

\end{prop}

The rest of this paper examines   certain geometric properties of the very particular case of 
$|St(a)| \leq d$, that will be used in \S \ref{sec_motive} to   prove that $\cM^{St(a)}$ define 
objects in the category of mixed Tate motives in these cases.

\subsection{Forgetful morphism} 
\label{sec_forget}

Consider the    trivial fibration   
\begin{eqnarray}  \label{eqn_forget1}
\begin{diagram}
(\A^d \times \A^d)^{Ver(\kappa)} \times \A^d	&    \rTo^{\psi_{b}}	 	&	 (\A^d \times \A^d)^{Ver(\kappa) \setminus \{b\}} \times \A^d  \\
(  {\bf z}^k  ,  {\bf w}^k  \mid k \in Ver(\kappa)) \times  {\bf w}^a	  & \rMapsto 
			&   (  {\bf z}^k  ,  {\bf w}^k  \mid k \in  Ver(\kappa) \setminus \{b\} ) \times  {\bf w}^a
\end{diagram}
\end{eqnarray}
whose fibres are $\A^d \times \A^d = \{ (  {\bf z}^b  ,  {\bf w}^b  ) \}$. 
The restriction of  $\psi_b$ to the configuration space $ \cM^{St(a)} $ of hyperplane  arrangements in almost general position
provides us {\it  the forgetful morphism}  
\begin{equation}  \label{eqn_hyperplane_forget}
	\begin{array}{rcl} 
	\psi_b: \cM^{St(a)} 			& 	\to		& \cM^{St(a) \setminus \{f\} } \\	
	(   P_{e}  \mid e \in St(a))	&	\mapsto	& (  P_{e}  \mid e \in St(a) \setminus \{f\})
	\end{array}
\end{equation} 
where $\{f\} = St(a) \cap St(b)$, i.e., $\partial_\kappa(f) = (ab)$. The image  $\psi_b (  P_{e}  \mid e \in St(a))$ of the 
hyperplane arrangements  
$ (  P_{e}  \mid e \in St(a)) $ is obtained by forgetting the hyperplane $ {P}_{f}$. 

\subsubsection{Fibres of the forgetful morphism when $|St(a)| \leq d$}  
\label{sec_forget_fiber}

We decompose   the forgetful morphism by using the following projections
\begin{eqnarray}  
\begin{diagram} \label{eqn_composition}
	(\A^d \times \A^d)^{Ver(\kappa)} \times \A^d	& \rTo^{\rho_{b}} & 	
		 (\A^d \times \A^d)^{Ver(\kappa) \setminus \{b\}}  \times  \A^d \times \A^d 
			& \rTo^{\phi_b}	&	 (\A^d \times \A^d)^{Ver(\kappa) \setminus \{b\}} \times \A^d \\
	\uInto & & \uInto & &  \uInto \\		
	\cM^{St(a)} &\rTo & \cN^{st(a)}& \rTo & \cM^{St(a) \setminus \{f\} } 
\end{diagram}
\end{eqnarray}
Here, the morphism $\rho_{b}$ forgets ${\bf z}^b$, that is  in fact   the product 
$$
\rho_{b} :=   \pi_{\kappa} \times \text{id} :( (\A^d \times \A^d)^{Ver(\kappa)}) \times \A^d \to 
((\A^d \times \A^d)^{Ver(\kappa) \setminus \{b\}} \times \A^d ) \times \A^d
$$ 
of the projection defined in \eqref{eqn_full_forgetful-1} with the identity morphism  of the last  factor 
$\A^d = \{ {\bf w}^a  \}$. 
The morphism $\phi_b$ forgets ${\bf w}^b$. 
The restriction of the composition $ \phi_b \circ \rho_b$
onto the configuration space $\cM^{St(a) } \subset (\A^d \times \A^d)^{Ver(\kappa)} \times \A^d$ gives the forgetful morphism 
$\psi_b$ defined in \eqref{eqn_hyperplane_forget}.

\smallskip
\paragraph{\bf Step 1. Fibers of $\phi_b$.}
 
Let $|Ver(\kappa)| \leq d$.  Consider a trivial bundle 
$$
   \cM^{St(a) \setminus \{f\}}  \times \A^d \to \cM^{St(a) \setminus \{f\}}
$$ 
as the restriction of  $\phi_b$ in \eqref{eqn_composition}.

For any $I \subset St(a) \setminus \{f\}$, there is a rank-$|I|$ 
subbundle $\cW_I$ whose fibers over 
$ (  {\bf z}^k  ,  {\bf w}^k   \mid k \in St(a) \setminus \{f \} ) \times  {\bf w}^a \in \cM^{St(a) \setminus \{f\}}$
consists of the space $\A^{|I|}$ which is the   vector space spanned by   
$({\bf w}^a - {\bf w}^k)$ for $k \in I$, that are the normals to the hyperplanes $P_e$ in $\cA_{St(a) \setminus \{f\}}$.

\begin{lem} \label{lem_fiber1}
The bundle $\cW_I$ over the variety $\cM^{St(a) \setminus \{f\}}$ is trivial for all $I \subset St(a) \setminus \{f\}$. 
\end{lem}

\begin{proof}
Above description of $\cW_I$ means precisely that the map
\begin{eqnarray*}
 \cM^{St(a) \setminus \{f\}}  \times \K^{|I|} & \to &   \cW_I \\
( (  P_{e}  \mid e \in St(a) \setminus \{f\}) \ , \  (c_k \mid k \in I  )) & \mapsto  &
( (  P_{e}  \mid e \in St(a)  \setminus \{f\}) \ ,\  (\sum_{k \in I} c_k ({\bf w}^a - {\bf w}^k) )  )
\end{eqnarray*}
is an isomorphism and provides the trivialization that is needed. 
\end{proof}

 The bundle $\cW_I$ parameterizes the pairs 
 $(  \cA_{St(a) \setminus \{f\}} , {\bf w}^b - {\bf w}^a )$ where the normal ${\bf w}^b - {\bf w}^a$ of
 forgetten hyperplane $P_f$ lies in the vector space spanned by the normal vectors 
 $({\bf w}^a - {\bf w}^k)$ for $k \in I$.
Such $P_f$'s should be in the complement of the configuration space due to \eqref{eqn_generic}.
Therefore, we will be interested in its complement  below.

\smallskip
\paragraph{\bf Step 2. Fibers of $\rho_b$.} 
The condition \eqref{eqn_generic} implies that the image of $\cM^{St(a)}$ under
$\rho_b$ is contained in the complement 
\begin{equation} \label{eqn_inter}
\begin{array}{ccl}
\cN^{St(a)} & :=  &( \cM^{St(a) \setminus \{f\}}  \times \A^d )  \setminus   \cW^{St(a)}   \\
 &  = & \{ (  \cA_{St(a) \setminus \{f\}} , {\bf w}^b - {\bf w}^a ) \mid 
 		{\bf w}^b - {\bf w}^a \not\in \left\langle {\bf w}^k - {\bf w}^a  : k \in I \ \& \ |I| \leq d \right \rangle   \} \
 \end{array}
\end{equation}
of the union $\cW^{St(a)}: = \bigcup_I \cW_I$. However, it is not clear whether the
image of this morphism covers $\cN^{St(a)}$. In the  following, we observe that the fibres of
$\rho_i$ are complements of hyperplane arrangements hence non-empty.

\begin{figure}
\begin{center}
\includegraphics[scale=.9]{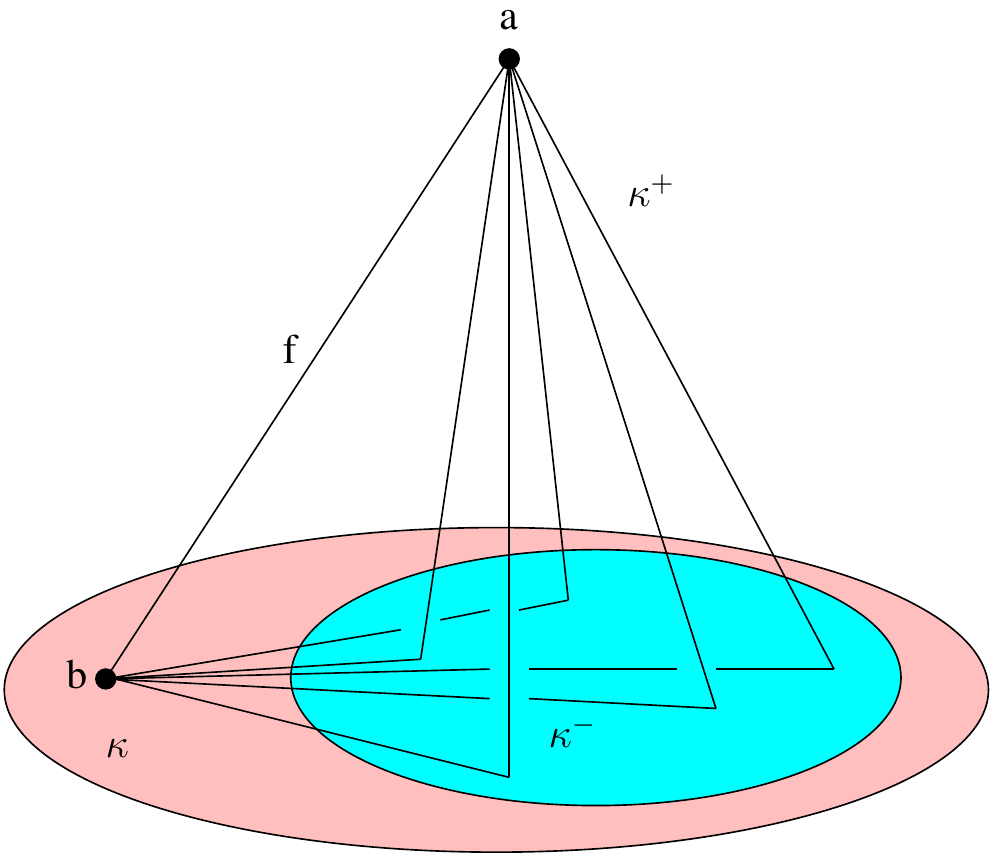}
\caption{ The complete graphs $\kappa^+$,   $\kappa$, and $\kappa^-$.  
\label{fig_forget}}
\end{center}
\end{figure}

\begin{lem} \label{lem_fiber2}
The fibre $\rho_i^{-1}({\bf n})$ over a point ${\bf n}: = (  \cA_{St(a) \setminus \{f\}} , {\bf w}^b - {\bf w}^a ) \in \cN^{St(a)}$
is the complement of a  hyperplane arrangement
$$
  \cB_{St(b)} = (   H_e \in \P^d \mid e \in ( St(b) \setminus \{f\} )   \subset Ver(\kappa))
$$    
where 
$$
 H_e := \{ {\bf z}^b  \mid q^{bj} = ({\bf z}^b - {\bf z}^j) \cdot ({\bf w}^b - {\bf w}^j )  = 0 \ \text{and} \ (bj) = \partial_\kappa(e)\}. 
$$
	Moreover, for any $I \subset St(b)$ (since we assumed  $|I| \leq |St(b)| \leq d$), 
	 the subarrangements $  \cB_I = (  H_{e}  \mid e \in I)$ are in general position, i.e., 
	\begin{equation*}
	  H_{I}  :=  \bigcap_{{e \in I }   }   H_{e} = \A^{d - |I|}.
	\end{equation*}

\end{lem}

\begin{proof}
This statement is just an iteration of the case (3) in  Lemma \ref{lem_full_fibre} where 
$\kappa, \kappa^-$ is a pair of complete graphs such that $\kappa^-$ is obtained from
  obtained from $\kappa$ by removing its vertex $b$  and its star $St(b)$, see Figure \ref{fig_forget}. 

All we need to make sure is that
  the elements   ${\bf n} \in\cN^{St(a)}$ is in the complement of the $  Z_{\kappa^-} \times \A^d$  
(which simply follows
from the fact that $\in\cN^{St(a)}$  fibres over $\cM^{St(a) \setminus \{f\}}$ which is in  the complement of $  Z_{\kappa^-} \times \A^d $) 
and $\bigcup_{e  \in St(b) } \{ {\bf w}^i  =  {\bf w}^a   \}$ (that is guaranteed as this locus lies inside  $\cW^{St(a)}$ (see,
Step 1 above)).
\end{proof}

 Consider the trivial bundle  $\cN^{St(a)} \times \A^d \to \cN^{St(a)}$ as the restriction of of  $\rho_b$ in 
 \eqref{eqn_composition}. For any $I \subset St(a) \setminus \{f\}$, there is a rank-$|I|$ 
subbundle $\cR_I$ whose fibers over 
$  (  \cA_{St(a) \setminus \{f\}} , {\bf w}^b - {\bf w}^a )$
is the intersection  $  R_{I}  :=  \bigcap_{{e \in I }   }   H_{e} = \A^{d - |I|}$ of the hyperplanes in Lemma \ref{lem_fiber2}.

 \begin{lem} \label{lem_trivial_again}
Let $\cR_I \to \cN^{St(a)}$ are trivial for all $I \subset St(a) \setminus \{f\}$.
 \end{lem}

 The proof this Lemma is same as Lemma \ref{lem_fiber1}.

\section{Motive of the Feynman quadric} 
\label{sec_motive}

In  this section,  we use the results from   \S \ref{sec_degeneration}
and \S \ref{sec_forget}
to  prove that the motive of  configuration space $\cM^{St(a)}$ is mixed Tate 
when $|St(a)| \leq d$.  We then 
prove that the Feynman quadrics  also give mixed Tate motives in corresponding cases
via  \S \ref{sec_projection_extended}.

\subsubsection{Assumptions for the induction} \label{sec_assumptions}

We prove our main theorem  by   induction on the number vertices of the complete graphs. 
We assume that the Feynman quadrics $Z_\kappa$ define mixed Tate motives  for all $|Ver(\kappa)| \leq n$.  
This assumption implies that the motives of configuration spaces $\cM^{St(a)}$ are also
mixed Tate   for all such $\kappa$ with  $a \in Ver(\kappa)$ (for details, see Case   (2) in the proof of 
Main Theorem). These  assumptions  are already verified for  the initial step of the induction, i.e., 1-edge
graph,  in \S \ref{sec_propagator}.

We will use the projections in \S \ref{sec_projection_extended} 
and \S \ref{sec_forget_fiber} to  examine the motive of the Feynman quadric $Z_{\kappa^+}$
for a bigger graph $\kappa^+$.

\subsection{The motive of the configuration space $\cM^{St(a)}$: The  case of $St(a) \leq d$} 

In this particular case, we will use the forgetful morphism \eqref{eqn_hyperplane_forget} and
its detailed analysis in \S \ref{sec_forget_fiber} to prove that
the configuration spaces define objects in the category of  mixed Tate motives.

In \S \ref{sec_forget_fiber}, we show that we need to consider the complement of 
certain stratified spaces to examine the motive of the configuration space $\cM^{St(a)}$.  
The following Lemma will be
useful for our purpose. 

\begin{lem} \label{lem_motiv_stratification}
Let $\{ X _i, i \in I \}$ be the set of irreducible components of a scheme $ X$ such that the
motives $\m^c( X_J)$ of the intersections 
$$
X_J : =  \bigcap_{j \in J } X_j
$$
are mixed Tate for all $J \subseteq I$. Then, the motive $\m^c(X)$ is also mixed Tate. 
\end{lem}

\begin{proof}
For inclusions 
$$
\bigcup_{ |J'| = k+1 } X_{J'} \hookrightarrow \bigcup_{|J| = k }  X_J, 
$$ 
we consider  the following canonical distinguished  triangle
\begin{eqnarray*}  
\begin{diagram}
\m^c(\bigcup_{ |J'| =k+1}  X_{J'} )   & \rTo &  \m^c(\bigcup_{ |J| = k}  X_{J}) & \rTo & 
		\m^c(\bigcup_{|J| = k} X_J \setminus \bigcup_{J' \supset J}  X_{J'} )
		& \rTo &   \m^c(\bigcup_{ |J'| = k+1  }  X_{J'} )[1].  
\end{diagram}
\end{eqnarray*}
If we assume that the motive $ \m^c(\bigcup_{|J| = k} X_J \setminus \bigcup_{J' \supset J}  X_{J'} )$ 
of the complement is mixed Tate for all  $k$, then we can claim that the motive of $X = \bigcup_i X_i$ is mixed 
Tate by induction on $k$: 
The case $k= |I|$ simply follows from the assumption that $\m^c(X_I)$ is mixed Tate.   As we assume that 
$ \m^c(\bigcup_{|J| = k} X_J \setminus \bigcup_{J' \supset J} X_{J'} )$ is mixed Tate, and know that
$\m^c(\bigcup_{ |J'| = k+1} X_{J'} )$ is mixed Tate from the previous step, the motive $\bigcup_{ |J| =k} X_{J}$
must be mixed Tate due to the above distinguished triangle. The final step $k=1$ gives us the motive 
$\m^c(\bigcup_i X_i)$ of $X$.

 It remains to prove our assumption, that $ \m^c(\bigcup_{|J| = k} X_J \setminus \bigcup_{J' \supset J} X_{J'} )$ 
 is mixed Tate for all  $k$. We simplify the problem by observing that 
$$
\m^c(\bigcup_{|J| = k} X_J \setminus \bigcup_{J' \supset J} X_{J'} ) = 
\bigoplus_{|J| = k} \m^c( X_J \setminus \bigcup_{J' \supset J} X_{J'} ).
$$
Notice that $\m^c(\bigcup_{J' \supset J} X_{J'})$ is mixed Tate due to induction hypothesis.
 Therefore, we can employ the   distinguished  triangle for the inclusion 
$\bigcup_{J' \supset J} X_{J'} \hookrightarrow X_J$
\begin{eqnarray*} 
\begin{diagram}
\m^c(\bigcup_{J' \supset J} X_{J'})    & \rTo &   \m^c( X_{J}) & \rTo &   
		\m^c( X_J \setminus \bigcup_{J' \supset J} X_{J'} )
		& \rTo &   \m^c(\bigcup_{ J'} X_{J'} )[1]. 
\end{diagram}
\end{eqnarray*}
to complete the proof. 
%
\end{proof}

\begin{lem} \label{lem_motiveN}
The motive of  $\cN^{St(a)}$ lies in the category of mixed Tate motives. 
\end{lem}

\begin{proof}
From  \S1.2.3 of \cite{Bloch} we know that the motive $\m(X \times \A^k) = \m(X)[-k] (2k)$ is obtained from  
$\m(X)$ by Tate twists and shifts, hence the bundles $\cW_I$ discussed in Lemma \ref{lem_fiber1} are in the
subcategory of mixed Tate motives inside the Voevodsky's category. Lemma \ref{lem_motiv_stratification}
implies that their union $\cW^{St(a)}$ also must be a mixed Tate motives. Therefore, we conclude that its 
complement $\cN^{St(a)}$ inside a mixed Tate motive $\cM^{St(a) \setminus \{f\} } \times \A^d$  (due to 
the induction assumption) is also a mixed Tate
motive. 
\end{proof}

The exact same line of arguments applies to the motive of the configuration space $\cM^{St(a)}$:

\begin{prop} 
\label{prop_deepest1}
If  $|St(a)|  \leq d$, then the configuration space $\cM^{St(a)}$ defines an object in the category  of 
mixed Tate motives.
\end{prop}

\begin{proof}
Due to Lemma \ref{lem_fiber2}, the configuration space $\cM^{St(a)}$ is the complement of 
the  union $\bigcup \cR_I$ inside $\cN^{St(a)} \times \A^d$.   According to Lemmata 
\ref{lem_trivial_again} and \ref{lem_motiveN},  the strata  $\cR_I$ 
are trivial $\A^{|I|}$-fibrations over a mixed Tate base $\cN^{St(a)}$.   
Hence, each  stratum $\cR_I$ defines a mixed Tate motive.  Therefore, the motive of  their union 
$\bigcup \cH_I$ is also mixed Tate due to Lemma \ref{lem_motiv_stratification}. 
The ambient space $\cN^{St(a)} \times \P^d$ is mixed Tate due to Lemma \ref{lem_motiveN},
so is its complement $\cM^{St(a)}$. 
\end{proof}

 \subsection{The motive of the Feynman quadric}
 Let $\kappa^+$ be a complete graph with $|Ver(\kappa)|  = n+1  \leq d+1$.

 \smallskip
 \noindent 
 {\bf Main Theorem:}  {\it   
The Feynman quadric  
$  Z_{\kappa^+}$ defines an object in the category of mixed Tate motives. }

\begin{proof}
We use  the projection
 \begin{equation*} 
	\begin{array}{rcl}	
	 \pi_{\kappa^+}:    Z_{\kappa^+}	& \to		& (\A^d \times \A^d)^{Ver(\kappa)} \times \A^d   \
	\end{array}
\end{equation*}
which was studied in Lemma \ref{lem_full_fibre}. Remember  that we assume in  \S \ref{sec_assumptions} 
that $\m^c(    Z_{\kappa} )$ is mixed Tate motives  to implement the induction  on the number of vertices.

According to Lemma \ref{lem_full_fibre}, the Feynman quadric $  Z_{\kappa^+}$ is a union
of three pairwise disjoint  pieces:
\begin{enumerate}
\item A trivial fibration $\cA \to   Z_{\kappa} \times \A^d $ with fibres $\A^d$, 
\item A trivial  fibration $\cB \to \bigcup_{e_i \in St(a)} \{ {\bf w}^a  = {\bf w}^i  \mid \partial_{\kappa^+}(e_i) = (ai) \}$ with fibres $\A^d$, 
\item The ``universal"  family over the configuration space  $\cM^{St(a)}$ of hyperplane arrangements.
\end{enumerate}

Case (1). It can be shown that the motive $\m^c(\cA)$ is mixed   Tate via  \S1.2.3 of \cite{Bloch}. 

Case (2). To be able use  the same arguments for trivial fibrations, one needs to check the motive of the mutual intersections
of the diagonals $\{{\bf w}^a  = {\bf w}^i    \} $ and their intersections with $  Z_\kappa \times \A^d$ in 
$(\A^d \times \A^d)^{Ver(\kappa)} \times \A^d$.  However, the intersections 
$\{{\bf w}^a  = {\bf w}^i   \}  \cap \{ {\bf w}^a  = {\bf w}^j   \} = \{{\bf w}^a  = {\bf w}^i    = {\bf w}^i  \}  \subset \cH^{ij} \times A^d$
are contained in  $  Z_\kappa \times \A^d$. Therefore, we only need to the check the motive of
pairwise disjoint subspaces $\{{\bf w}^a   = {\bf w}^i   \} \setminus ((  Z_\kappa  \times \A^d) \cap \{ {\bf w}^a  =  {\bf w}^i    \} )$.  

Note  that 
 $(   Z_\kappa  \times \A^d) \cap    \{   {\bf w}^a  =  {\bf w}^i  \}$ is isomorphic to $  Z_\kappa$
as this intersection as be given as the graph of the map
\begin{eqnarray*}  
\begin{diagram}
  Z_\kappa 	&	 \rInto	& 	   Z_\kappa  \times \A^d   \\
( {\bf z}^m  ;   {\bf w}^m   \mid m \in Ver(\kappa))  &  \rMapsto &   ( {\bf z}^m ;  {\bf w}^m  \mid m \in Ver(\kappa)) \times   {\bf w}^i,
\end{diagram}
\end{eqnarray*}
which shows that $\m^c((  Z_\kappa   \times \A^d) \cap    \{   {\bf w}^a  =  {\bf w}^i  \} )$ and, therefore,
the motives $\m^c(\{ {\bf w}^a   = {\bf w}^i    \} \setminus ((  Z_\kappa   \times \A^d) \cap \{ {\bf w}^a   = {\bf w}^i   \} ))$ 
of the compliments are mixed Tate motives, so are the locally trivial $\A^d$-bundles $\cB$ over them.

Case(3).  The remaining part of $  Z_{\kappa^+}$ after Case (1) and (2) is the subspace 
$  \pi_{\kappa^+}^{-1} (\cM^{St(a)}) \cap   Z_{\kappa^+}$.  
This subspace is stratified by the trivial fibration  $\cR_{I}$ as in Lemma  
\ref{lem_trivial_again}.  Following the same steps in Proposition \ref{prop_deepest1}, we show that 
the remaining part, the universal family $  \pi_{\kappa^+}^{-1} (\cM^{St(a)}) \cap   Z_{\kappa^+}$ over 
the configuration space $\cM^{St(a)}$   also defines
an object in the category of mixed Tate motives.

Finally,  we  patch these three pieces together via the distinguished triangles and show that the Feynman quadrics  
$\bigcup_{e  \in Edg(\kappa^+)}   \cH^{e}$ indeed defines mixed Tate motives as the same was true
for one-edge graph, see Proposition \ref{prop_motive1} .
\end{proof}




 \section{Corollaries and   ramifications}
 \label{sec_corollaries}

 \subsection{Graph hypersurfaces vs  Feynman quadrics} 
 The {\it graph polynomial }of $\Gamma$ is given as 
 $$
 \Psi_\Gamma : = \sum_T \prod_{e \not\in T} \alpha_e
 $$
where $\alpha_e$'s are the variables associated to the edges of $\Gamma$
and the summation runs over the set of spanning trees $T \subset \Gamma$. 
We note that the graph hypersurfaces $X_\Gamma = \{ \Psi_\Gamma  =0\} \subset \A^{Edg(\Gamma)}$
do not depend on the dimension of the spacetime, hence the dimension of
the corresponding momentum space. On the other hand, we  can choose the 
dimension high enough so that the    residues of the form  \eqref{eqn_integral_complete} 
for can be interpreted in terms of the periods of a mixed Tate motive, the motive of
the Feynman quadric $Z_{\kappa_n}$.

 The contrast between our  main results  and the results on graph hypersurface 
 (see \cite{BD, BS, D} for counter examples and \cite{Mar} for a general detailed 
 account on the subject) draws our attention 
to two main unknown knowns:

\subsubsection{Transferring between the position and the momentum space is not as direct as hoped } 

The Fourier transform provides an isomorphism between the state spaces (i.e., the spaces of  square
integrable  functions) in the momentum and position spaces. However, this isomorphism utilizes  the cut-off function, that is likely to prohibit us
from transferring the algebraic therefore motivic structures directly.

\subsubsection{The dependence on regularization scheme}
The obstacle due to (the transcendental nature of) the Fourier transform however can be overcomed by 
directly considering the Feynman quadric in momentum space setting. 
Our recent computation in  \cite{Cey} shows that momentum space Feynman quadrics also define
mixed Tate motives. This concludes  that the motives arising in momentum space depends on
the choice of the regularization procedures. 

 The parametric formulation of Feynman integrals uses an integral presentation of propagators and
 the Fubini theorem to change the orders of (divergent) integrals.   
 These results indicate that the integral 
 presentation of propagator and introducing Schwinger parameters are transcendental in some sense
 as they replace mixed Tate motives, the Feynman quadrics, with non-mixed Tate ones, the graph
 hypersurfaces. 

The potential problems of  changing the orders of divergent integrals is implicitly noted in \cite{BEK}.
Bloch, Esnault and Kreimer considered the  locus of the  divergence of the momentum space propagator
as the zeros of the smooth quadric
$$
\{ x_1^2 + \cdots + x_d^2 = 0 \} \subset \A^d. 
$$
The cohomology of the union of these momentum space quadrics as well as their periods have been 
consider in  \S 10 of \cite{BEK} in case of logarithmically divergent primitive graphs, i.e., where the period 
integral is convergent and Schwinger trick cannot cause above mentioned problems.

In addition to these, there is a recent observation on mixed Tate motives which may also play
a role in Feynman integrals: 

\subsubsection{Tate motives are more elusive than one anticipates}
A very recent paper \cite{Borisov} observes an unexpected property of the mixed Tate motives: 
The class of the affine line is a zero divisor in the Grothendieck ring of algebraic varieties. 
This result is counter intuitive and indicates  that  the mixed Tate motives are in fact more 
sophisticated than they are generally  depicted.

 \subsection{Feynman quadrics for general Feynman graphs}
 The main construction of this paper fails at a number places for graphs other than the complete graphs: 
 Firstly, if we remove the restriction $|Ver(\Gamma)| \leq d+1$, we need to consider 
 hyperplane arrangements in almost general position, see \S \ref{sec_degeneration}. 
 For more general Feynman graphs, the    arrangements having parallel hyperplanes appear 
 in the picture as the corresponding case in Lemma \ref{lem_full_fibre} uses
 the fact that each pair of vertices connected by an edge to establish the  non-parallelity. 
 The same statements  are  also true for graph with external structures. 
 
These cases require more elaborate study of the motives of  the stratification of the respective 
configuration spaces. We plan to write these details in a separate paper. The result remain the 
same however the motives of the Feynman quadrics in such cases  fail to be unramified
 over $Spec(\Z)$ due to trivial reasons: The arrangements containing more than
 $(d+1)$ hyperplanes almost never unramified, and we need to consider the spaces of
 them.

 \subsection{The signature doesn't change anything}
The constructions in this paper focuses on Euclidean case, the metric with signature
$(+, +, \dots , +)$.  However, if we change the signature, we 
only need to replace our simple Feynman quadric \eqref{eqn_quadric} with 
$$
q'({\bf z}, {\bf w} ) = z_1 \cdot w_1 + \cdots + z_k \cdot w_k - z_{k+1} \cdot w_{k+1}
-  \cdots  - z_d \cdot w_d = 0.
$$
The rest of the statements and the proofs remain the same.

 \subsection{Other geometries and motives for position space Feynman integrals}

While the main problem seems well defined  at  the first glance, it hides an ambiguity between the lines.
As a period, the integral \eqref{eqn_feynman_int} must be an integral over a semialgebraic set. However, as
defined in  \eqref{eqn_feynman_form}, there is no canonical choice among the complexifications. In fact,  there any infinitely
many complexifications of the spacetime $X$ such that restriction of the complexified propagator $G({\bf x})$ 
to its real part would provide the propagator on $X$.

As a simplest example, one can start with any given complexification  of the spacetime $X$
and blow it up along subvarities having no real points. Such blowups do not change the real locus of
$X^{|Ver(\Gamma)|}$ and the propagator, therefore, provide different complexifications.  
However, there are more interesting cases in which the complexification remains the same.  Here is a sample 
of real structures on $\A^{2d}(\C)$ that each lead to quite different loci of divergences 
despite that their restriction to real locus $\A^{2d}(\R)$ are the same. 
\begin{eqnarray*}
\begin{array}{cccc}
 &  \text{Real structure on } 	\A^{2d}(\C)	  & \text{Propagator}  & \text{Locus of divergence} \\
 (I) & (z_1,\dots,z_{2d}) \mapsto (\bar z_1,\dots, \bar z_{2d})        
	&  \frac{1}{( \sum_{i=1}^d \epsilon_i (z_i)^2)^{d-1}} & \left\{  \sum_{i=1}^d \epsilon_i (z_i)^2   = 0 \right\}, \\
 (II) &(z_1,\dots,z_{2d}) \mapsto (\bar z_1,\dots, \bar z_{2d})     
	&  \frac{1}{(\sum_{i=1}^d |z_i|^2)^{d-1} } & \left\{  z_1 =  \cdots  = z_{2d} = 0 \right\}, \\	
 (III) & (z_1,\dots,z_{2d}) \mapsto (\bar z_{d+1},\dots,  \bar z_{2d}, \bar z_1, \dots, \bar z_d)        
	&  \frac{1}{ ( \sum_{i=1}^d  \epsilon_i  (z_{i} \cdot z_{d+i}) )^{d-1}} & \left\{  \sum_{i=1}^d  \epsilon_i  (z_{i} \cdot z_{d+i}) = 0 \right\}. \\
\end{array}
\end{eqnarray*}
where $(\epsilon_i = \pm 1)$ is the signature of the metric. 

The loci of divergences are all algebraic varieties  each having distinct geometric properties: While the loci of divergence in  (I) and (II) 
are smooth, the locus of divergence in (III) has a singularity at the origin. In the case (I) and (III), they are hypersurfaces however 
it is a codimension $d$ subvariety in the case of (II). Finally, being mixed Tate motive is by definition in the case of 
(II), the same property   is quite nontrivial in the case (III) as we have seen in this paper.

 \subsubsection{Feynman integrals as configuration space integrals}
 
 In a series of paper \cite{CeyMar1, CeyMar2, CeyMar3}, the  Feynman integrals are studied as the homological 
 pairings in the configuration spaces of points defined via the propagator in  (II). It is quite easy to observe that these configuration 
 spaces of points can be nicely stratified and desingularized, and they all define  mixed Tate motives. However, the
 propagator defined in (II)  is a real valued distribution and is not algebraic, i.e., the corresponding Feynman integrals 
 cannot manifest themselves  as periods directly.

\subsection*{Acknowledgments} 
I would like  to thank to Matilde Marcolli for uncountable   many  discussions over the years. 
I also wish to thank to Paolo Aluffi for pointing out the ambiguities in the earlier version of this
paper and
helping on clarifying the statements, to G\"unther Harder for  discussion on ramifications 
of motives, and, to Yuri Manin for many inspiring discussions in particular for attracting my
attention to recent developments \cite{Borisov}. 
Thanks are also due to
K\"ur\c{s}at Aker and Atabey Kaygun for their help and suggestions while drafting this paper. 

Parts of this work have been carried out during visits at Max-Planck Institut f\"ur Mathematik in Bonn. 
The   author is partially supported by Marie-Curie Career Integration Grant  PCIG11-GA-2012-322154 and 
FNR Open Grant QUANTMOD O13/5707106.

\end{document}